\RequirePackage[l2tabu, orthodox]{nag}

\documentclass[a4paper,12pt,french,english,leqno]{article}

\usepackage[utf8x]{inputenc}

\usepackage[T1]{fontenc}
\usepackage{babel}
\usepackage[margin=2cm]{geometry}

\usepackage{amsmath,amssymb,amsthm,amscd}
\usepackage{euscript,mathrsfs}
\usepackage{mathtools}
\usepackage[shortlabels]{enumitem}
\usepackage{array}

\usepackage[svgnames,x11names]{xcolor}
\usepackage[colorlinks=true]{hyperref}
\usepackage{txfonts}
\usepackage{pgf}
\usepackage{graphicx}

\def\ge{\geqslant}
\def\le{\leqslant}

\def\N{\mathbb{N}}

\def\C{\mathbb{C}}

\def\setminus{\smallsetminus}
\def\emptyset{\varnothing}

\def\A{\mathbb{A}}
\def\angle#1{\langle#1\rangle}
\def\CC{\mathcal{C}}
\def\Cred{C_{\mathrm{red}}}

\def\iso{\simeq}
\def\L{\mathcal{L}}
\def\m{\mathfrak{m}}
\def\mapto{\longrightarrow}
\def\maptodots{\mathrel{{\cdot}{\cdot}{\cdot}{\rightarrow}}}
\def\maptoinj{\hookrightarrow}
\def\maptosurj{\twoheadrightarrow}
\def\O{\mathcal{O}}
\def\P{\mathfrak{P}}
\def\PP{\mathbb{P}}
\def\tC{\tilde C}
\def\tH{\tilde H}

\DeclareMathOperator\height{ht}
\DeclareMathOperator\ord{ord}
\DeclareMathOperator\QF{QF}
\DeclareMathOperator\Spec{Spec}

\newtheorem{lemma}{Lemme}[section]
\newtheorem{definition}[lemma]{Définition}

\newtheorem{coro}[lemma]{Corollaire}
\newtheorem{propo}[lemma]{Proposition}

\newtheorem{prf}{\it{Preuve}}

\newenvironment{demo}{\begin{prf}\rm}{\hfill$\Box$\end{prf}}
\def\cad/{c'est-à-dire,}
\def\ie/{c'est-à-dire,}
\def\pb/{point base}
\def\pbs/{points base}

\def\English{\foreignlanguage{english}}

\def\frac#1#2{#1/#2}

\title{Existence des diviseurs dicritiques,\\d'après S.S. Abhyankar.}
\author{Vincent Cossart, Mickaël Matusinski, Guillermo Moreno-Socías}
\date{1\textsuperscript{er} mars 2013}

\begin{document}
\maketitle

\begin{flushright}
\textit{À la mémoire de Shreeram Shankar Abhyankar (1930--2012)}
\end{flushright}

\bigskip

\begin{center}
\bfseries\large
INTRODUCTION
\end{center}

Soit $(f,g)$ un couple de polynômes de $\C[X,Y]$.
On considère leurs homogénéisés $(F,G)\in \C[X,Y,Z]$ premiers entre eux et de même degré.

On a alors une fonction de 
$\PP^2_{\C} \mapto \PP^1_{\C}$, $(X;Y;Z)\mapto (F(X,Y,Z);G(X,Y,Z))$. 
Cette fonction n'est bien sûr pas définie aux \pbs/ du pinceau $\angle{F,G}$, 
mais on peut la définir sur une surface obtenue à partir de $\PP^2_{\C}$ en éclatant les \pbs/. 
Les diviseurs dicritiques de $\angle{F,G}$ sont les diviseurs exceptionnels 
tels que l'application restreinte à ces diviseurs est surjective (voir Définition~\ref{defi=dicrit1}). 
Ces diviseurs ont un rôle crucial dans le problème jacobien~\cite{le-weber_jacobian-dim2}. 

On retrouve les diviseurs dicritiques chez d'autres auteurs sous des qualificatifs différents~: 
diviseurs horizontaux chez Campillo--Reguera--Piltant \cite[Definition 4]{campillo-piltant-reguera=cones-curves2} 
et, dans un cadre plus général, 
diviseurs associés à des valuations de Rees de l'idéal des \pbs/ du pinceau $\angle{F,G}$ 
chez I.~Swanson \cite[Definition 1.1]{swanson=rees-val}\cite[Ex. 14.18]{huneke-swanson=integral-closure}, 
voir~\ref{propo=Rees} ci-dessous. 


Abhyankar a donné une définition des diviseurs dicritiques qui généralise et algébrise la définition géométrique précédente 
dans le cas local (\cite[Note (5.6)]{abh=inversion-invariance-char-terms} et Définition~\ref{defi=dicrit} ci-dessous) 
et dans le cas polynomial (\cite[Definition (5.1)]{abh=inversion-invariance-char-terms} et Définition~\ref{defi=dicrit2} ci-dessous). 
En suivant son exposé \cite[Section 5]{abh=inversion-invariance-char-terms}, 
nous donnons des interprétations géométriques des diviseurs dicritiques et des preuves nouvelles de leur existence.

Nous remercions Olivier Piltant pour ses explications et ses nombreux croquis qui nous ont permis de donner une nouvelle généralisation~\ref{propo=Abhyankar-Luengo} du théorème d'Abhyankar--Luengo 
\cite[Theorems (1.1), (7.1), (7.2), (7.3)]{abh-luengo_dicrit-div} 
et une preuve géométrique de l'existence des dicritiques~\ref{propo=Rees}.

C'est un article de mise au point avec un point de vue résolument géométrique. 
Le seul résultat nouveau est~\ref{coro=Abhyankar-Luengo-geom} 
qui donne un éclairage géométrique au théorème d'Abhyankar--Luengo 
et généralise le théorème de connexité de \cite{le-weber_jacobian-dim2} p.~377.

\section{Cas local}

Tout au long de cette section, on note $R$ un anneau local régulier de dimension~$2$, 
$\m$ son idéal maximal et $K:=R/\m$ son corps résiduel. 
L'anneau de valuation $V$ désigne un diviseur premier de $R$, 
\ie/ un anneau de valuation discrète dominant $R$ avec extension résiduelle transcendante. 
On note $\m_v$ son idéal maximal et $K_v:=V/\m_v$ son corps résiduel. 
La projection canonique est $H_v:V\mapto K_v$, où $K$ est identifié à $H_v(R)$. 
Sous ces conditions, nous écrivons $K_v=K'(t)$ où $K'$ est la clôture algébrique relative de $K$ dans $K_v$, 
et $t$ est transcendant sur $K$. Et $\QF(R)$ désigne le corps de fractions de~$R$.
\begin{definition}\label{defi=dicrit}
Soit $z\in \QF(R)$, $z\ne 0$. 
On appelle \textbf{diviseur dicritique} de $z$ dans $R$ 
tout diviseur premier $V$ de $R$ tel que $z\in V$ et $H_v(z)$ est transcendant sur~$K$. 
\end{definition}

\begin{propo}\label{propo=finitude-dicrit}
Tout $z\in \QF(R)$ non nul a un nombre fini de diviseurs dicritiques, 
nombre qui est nul si et seulement si $z\in R$ ou $\frac{1}{z}\in R$. 
\end{propo}

\begin{demo}
Si $z\in R$ ou $\frac{1}{z}\in R$, alors $H_v(z)\in K$ pour tout diviseur premier $V$. 
Donc $z$ ne peut avoir de diviseur dicritique. 

Désormais, on suppose $z\notin R$ et $\frac{1}{z}\notin R$. 
Comme $z\notin R$ et $\frac{1}{z}\notin R$, on a $z=\frac{f}{g}$, 
fraction irréductible avec $f,g\in \m\setminus\{0\}$. 
On définit la suite d'éclatements suivante~:
$$W_0=\Spec(R) \longleftarrow W_1 \longleftarrow W_2 \longleftarrow ~~\cdots~~ \longleftarrow W_n$$ 
centrés en $x_i\in W_i$, $0\le i \le n-1$, 
$x_0=x=\m$, 
$x_i$ se projetant sur $x_{i-1}$, $0\le i \le n-1$
et tel que $I_i$ le transformé faible de $I_0=(f,g)$ ne soit pas principal en $x_i$ 
mais que $I_n$ soit principal en tout point de $X_n$ se projetant sur $x_{n-1}$. 
Montrons que le diviseur exceptionnel $E$ de $W_{n-1} \longleftarrow W_n$ est dicritique. 
En $x_{n-1}$, on note $I_{n-1}=(f_{n-1},g_{n-1})$ 
où 
$f_{n-1}=\frac{f}{M_{n-1}}$,
$g_{n-1}=\frac{g}{M_{n-1}}$ et
$M_{n-1}$ est un monôme de composantes exceptionnelles. 
On montre par récurrence sur $n$ que $f_{n-1}$ et $g_{n-1}$ 
sont premiers entre eux en tout point de $E$, avec $E\subseteq W_n$. 
Montrons que $f_{n-1}$ et $g_{n-1}$ sont de même ordre $\m_{n-1}$-adique. 
Sinon, par exemple $\ord_{\m_{n-1}}(f_{n-1})< \ord_{\m_{n-1}}(g_{n-1})$. 
On note $i=\ord_{\m_{n-1}}(g_{n-1})- \ord_{\m_{n-1}}(f_{n-1})$. 
Alors, en $x_n\in X_n$ sur le transformé strict $f_n$ de $f_{n-1}$, 
on a $I_n=(f_n,t^{i} g_n')$ où $t$ est une équation locale de $E$ 
et $g_n'$ le transformé strict de $g_{n-1}$. 
Si $I_n$ était principal en $x_n$, alors $f_n$ diviserait $g_n'$, 
et cela contredirait l'hypothèse que $f_{n-1}$ et $g_{n-1}$ sont premiers entre eux en $x_{n-1}$. 
Donc, en tout point de $E$, on a $I_n=(f_n,g_n)$ 
où $f_n$ et $g_n$ sont les transformés stricts de $f_{n-1}$ et $g_{n-1}$. 
Soit $y_1\in E$ avec $f_n(y_1)=0$ et $y_2\in E$ avec $g_n(y_2)=0$. 
Remarquons que $y_1\ne y_2$. 
Sinon, en $y_1$, $I_n$ étant principal, par exemple $f_n$ diviserait $g_n$; alors
$f_n$ et $g_n$ auraient une composante commune dans $\O_{W_n,y_1}$ 
et donc $f_{n-1}$ et $g_{n-1}$ en auraient une dans $\O_{W_{n-1},x_{n-1}}$.
On a $z(y_1)=0$ et $z(y_2)=\infty$; donc $E$ est dicritique.

Malheureusement, ce procédé ne donne pas tous les diviseurs dicritiques. 
Montrons néanmoins qu'il n'y en a qu'un nombre fini. 
Soit 
\begin{equation}\label{eq=1}
W_0=\Spec(R) \leftarrow\cdots\leftarrow W,
\end{equation}
où $W_{i-1}\leftarrow W_i$ est cette fois-ci 
l'éclatement de \emph{tous} les points fermés $y\in W_{i-1}$ 
où le transformé faible de $I$ n'est pas $\O_{W_{i-1},y}$. 
Il est connu que cet algorithme est fini \cite[Lemme 2.1.1.]{giraud=forme-normale-surf}.

Soit $V$ un diviseur dicritique pour $z$. 
Par le critère valuatif de propreté, $V$ a un centre $A$ sur $W$, 
\cad/ $R=\O_{W_0,x_0}\subseteq \O_{W,A}\subseteq V$ et $\m_v\cap \O_{W,A}=\m_A$. 
Si $A$ est le point générique d'une courbe~$E$, 
alors $E$ est exceptionnelle et $\O_{W,A}=V$. 
Sinon, $A$ est un point fermé et $z=\frac{f'}{g'}$ avec $(f',g')$ transformé faible de $(f,g)$ en $A$. 
Mais par construction de $W$, $f'$ ou $g'$ est inversible en $A$. 
Donc $z$ ou $\frac{1}{z}$ appartient à $\O_{W,A}$ : 
$H_v(z)$ ou $H_v(\frac{1}{z})$ appartient à $\O_{W,A}/\m_A$ qui est algébrique su $K=R/\m$, 
$V$ n'est pas dicritique pour $z$.

Les diviseurs dicritiques sont donc parmi les $\O_{W,\eta}$ 
avec $\eta$ point générique d'une composante exceptionnelle : il y en a un nombre fini.


Dans \cite{abh=inversion-invariance-char-terms}, 
on trouvera une autre preuve de~\ref{propo=finitude-dicrit} 
avec des arguments plus algébriques et ex\-trême\-ment informatifs d'Abhyankar. 
La rédaction étant très concise, nous proposons ici une nouvelle rédaction plus détaillée.

Montrons que :
\begin{lemma}%
\footnote{Nous remercions le rapporteur qui remarqua que l'hypothèse $R$ factoriel est suffisante pour ce lemme.}
Soient $R$ un anneau factoriel et $z\in \QF(R)$ avec $z\notin R$ et $1/z \notin R$.
$V$ est un diviseur premier de $R$.
Et étant donnée une variable abstraite $Z$, il existe un homomorphisme surjectif 
$$ h:R[z]\mapto K[Z]$$
défini par $z\mapsto Z$ et $x\mapsto H_v(x)$ pour tout $x\in R$. 
\end{lemma}
\begin{demo}
On reprend les notations de ci-dessus~: 
$z\notin R$ et $\frac{1}{z}\notin R$, 
on a $z=\frac{f}{g}$, fraction irréductible avec $f,g\in \m\setminus\{0\}$. 
Soit $\pi: R[Z] \mapto R[z]$ la surjection naturelle. 
Montrons que
$\ker(\pi)={\angle{f-gZ}}$. 
Bien sûr, $\ker(\pi)\supseteq {\angle{f-gZ}}$.
Soit $P(Z):=a_d Z^d+a_{d-1}Z^{d-1}+\dots+a_1 Z+a_0 \in \ker(\pi)$ non nul.
Il est clair que $d\ge 1$. 
Montrons par récurrence sur $d$ que $P \in {\angle{f-gZ}}$.
On a $a_d f^d +g (a_{d-1}f^{d-1}+\dots+a_0 g^{d-1})=0$, donc $g$ divise $a_d$, on a $a_d=gb$.
Puis $P(Z)=bZ^{d-1}(gZ-f)+Q(Z)$ avec $Q(Z)\in \ker(\pi)$, $\deg Q\le d-1$, 
donc $Q$ est dans $\angle{f-gZ}$ et $P$ aussi.

Les homomorphismes suivants :
\begin{alignat*}{2}
\tH_v &: R[Z] &{}\mapto{} & K[Z]\\
\pi &:R[Z] &{}\mapto{} &R[Z]/\angle{f-gZ}\iso R[z] 
\end{alignat*}
sont surjectifs et tels que $\angle{f-gZ} = \ker(\pi) \subseteq \ker(\tilde{H}_v) $. 

\begin{center}
\begin{pgfpicture}{-3.5cm}{-2.5cm}{3.5cm}{1cm}%
\pgfsetroundjoin%
\pgfsetstrokecolor{black}
\pgfsetlinewidth{0.2pt} 
\pgfsetarrows{-to}
\pgfxyline(-2,-0)(2,0)
\pgfsetarrows{-}
\pgfputat{\pgfxy(-2.8,-0)}{\pgftext{\color{black}\small $R[Z]$}}\pgfstroke
\pgfputat{\pgfxy(2.8,0)}{\pgftext{\color{black}\small $K[Z]$}}\pgfstroke
\pgfputat{\pgfxy(0,0.3)}{\pgftext{\color{black}\small $\tilde{H}_v$}}\pgfstroke
\pgfputat{\pgfxy(-0,-1.8)}{\pgftext{\color{black}\small $R[z]=\frac{R[Z]}{(f-gZ)}$}}\pgfstroke
\pgfsetarrows{-to}
\pgfxyline(-2.5,-0.3)(-0.5,-1.3)
\pgfsetarrows{-}
\pgfsetdash{{5pt}{3pt}}{0pt}
\pgfsetarrows{-to}
\pgfxyline(0.5,-1.3)(2.5,-0.3)
\pgfsetarrows{-}
\pgfputat{\pgfxy(-2.1,-1)}{\pgftext{\color{black}\small $\pi$}}\pgfstroke
\pgfputat{\pgfxy(2.1,-1)}{\pgftext{\color{black}\small $h$}}\pgfstroke
\end{pgfpicture}%

\end{center}

L'application $h$ est l'unique homomorphisme tel que $\tH_v=\pi\circ h$.
\end{demo}
Le lemme suivant nous permet d'affirmer que, 
dans le cas où $V$ est un diviseur dicritique de $z$ dans $R$ anneau local régulier de dimension~$2$, 
on a $\m[z]=\ker(h)$ premier dans $R[z]$, 
$\m[z]=\m_v \cap R[z]$ et
$\m[z]\cap R=\m$. 
De plus, si l'on note $S:=R[z]_{\m[z]}$ le localisé de $R[z]$ en $\m[z]$, 
alors $\dim(S)=1$.

\begin{lemma}
Avec les hypothèses et notations ci-dessus, on pose 
$\P := \m_v \cap R[z]$.
Soit $z \in V$ de valuation nulle.
\begin{enumerate}[(i)]
\item\label{lemme=i}
Si $H_v(z)$ est algébrique sur $K$,
alors
$\m[z] \subsetneq \P$, $\P$ est maximal et 
$\frac{R[z]}{\P}$ est une extension finie de~$K$.
\item\label{lemme=ii}
Si $H_v(z)$ est transcendant sur $K$, 
alors on a que $\P = \m[z]$, $\frac{R[z]}{\P} \iso K[Z]$, 
$\dim(R[z]) =2$ et $\dim(R[z]_{\P}) =\nobreak1$.
\end{enumerate}
\end{lemma}
\begin{demo}
Bien sûr, $\m[z] \subseteq \P$. 
On a une suite 
d'injections~:
$K \maptoinj \frac{R[z]}{\P} 
\maptoinj K_v$.
L'injection $K \maptoinj \frac{R[z]}{\P}$
est la composée des flèches naturelles 
$K \mapto K[Z]$ et
$K[Z] \maptosurj \frac{R[z]}{\P}$.
Le dernier morphisme est surjectif, nous le notons $s$. 

\ref{lemme=i}
$H_v(z)$ est algébrique sur $K$ dans $K_v$ 
si et seulement si la classe de $z$ dans $\frac{R[z]}{\P}$ est algébrique sur $K$ :
on a une relation 
$H_v(a_0)+H_v(a_1)H_v(z)\cdots+H_v(a_m)H_v(z)^m =0 \in K_v$ avec 
$a_i \in R$, $0 \le i \le m$, et $H_v(a_i)\ne 0$ pour au moins un $i$.
Il est clair que $H_v(a_0)+H_v(a_1)Z\cdots+H_v(a_m)Z^m \in \ker(s)\ne (0)$.
Le noyau de $s$ est un idéal maximal
et $\frac{R[z]}{\P}$ est une extension finie de $K$. 
Considérons les morphismes naturels $\pi: R[Z]\mapto R[z]$
et $\tH_v:R[Z]\mapto K[Z]$. On a :
\begin{equation*}
\begin{CD}
R@>>>R[Z]@>\pi>>R[z]@>>>V\\
@V{H_v}VV@V{\tH_v}VV@VVpV@VV{H_v}V\\
K@>>>K[Z]@>s>>R[z]/\P@>>>K_v
\end{CD}
\end{equation*}
avec 
$\pi^{-1}(\P)=\tH_v^{-1}(\ker(s))$ et
$\pi^{-1}(\m[z])=\m[Z]=\tH_v^{-1}(0)$. 
Comme $\tH_v$ est surjective et $\ker(s)\ne (0)$, on a $\m[z] \subsetneq \P$.

\ref{lemme=ii}
Supposons que $H(z)$ est transcendant sur $K$, 
ce qui équivaut à dire que la classe de $z$ dans $\frac{R[z]}{\P} $ est transcendante sur $K$.
Alors le noyau de $s: K[Z] \mapto \frac{R[z]}{\P}$ est $(0)$, et $s$ est un isomorphisme.

Pour conclure, soit $x:=a_0+a_1z\cdots+a_mz^m \in R[z]\setminus \m[z]$.
Il existe $i, 0 \le i \le m$, tel que $a_i\notin \m$.
Alors $H_v(a_0)+H_v(a_1)Z\cdots+H_v(a_m)Z^m \ne 0 \in K[Z]$. 
Puisque $s$ est un isomorphisme, on a 
$H_v(a_0)+H_v(a_1)H_v(z)\cdots+H_v(a_m)H_v(z)^m \ne 0 \in \frac{R[z]}{\P}$,
donc $x\notin \P $. 
Ainsi $\m[z]=\P$ et comme $K[Z] \iso \frac{R[z]}{\P}$, on a que $\P$ n'est pas maximal.

On a donc dans $R[z]$ la chaîne d'idéaux premiers~: 
$(0)\subsetneq \m[z]=\P \subsetneq M$, 
où $M$ est un idéal maximal. 
D'autre part,
$\frac{R[Z]}{\angle{f-gZ}}\iso R[z]$, 
donc $\dim(R[z])\le \dim(R)=2$. 
Donc $\dim(R[z])=2$ et $\dim(R[z]_{\P})=1$.
\end{demo}


Pour conclure la preuve de la Proposition~\ref{propo=finitude-dicrit}, 
$S:=R[z]_{\m[z]}$ est l'anneau de la courbe $\CC$ d'équation $(f-gZ)$ 
de l'anneau régulier $R[Z]$ de dimension~$2$. 
Autrement dit, $\CC$ est la courbe générique du pinceau $(f,g)$ 
au sens de \cite[p.517]{campillo-piltant-reguera=cones-curves2}.
Sa normalisation est de Krull, noethérienne, à fibres finies, 
d'après les Théorèmes 33.10 et 33.12 de \cite{nagata_local-rings}. 
Ainsi, si l'on note $T$ la clôture intégrale de $S$ dans $\QF(R)$, on a :
$$T=V_1\cap\cdots\cap V_e$$
où $e\in\N^*$ et les $V_i$ sont des anneaux de valuation discrète deux à deux distincts de $\QF(R)$. 
Ce sont précisément les diviseurs dicritiques de $z$ dans $R$.
\end{demo}

\section{Diviseurs dicritiques et valuations de Rees}

\begin{definition}
Soit $I$ un idéal de $R$. 
L'ensemble des \textbf{valuations de Rees} de $I$ 
est le plus petit ensemble $\{V_1,\ldots,V_e\}$ d'anneaux de valuation vérifiant :
\begin{enumerate}[1.]
\item
Les $V_i$ sont noethériens et ne sont pas des corps.
\item
Pour tout $i$, il existe un idéal premier minimal $\P_i$ de $R$ 
tel qu'on a la chaîne d'anneaux $\frac{R}{\P_i}\subseteq V_i\subseteq \QF\left(\frac{R}{\P_i}\right)$.
\item
Pour tout $n\in\N $, la clôture intégrale $\overline{I^n}=\bigcap_{i=1}^e(I^nV_i)\cap R$.\\
N.B. : il s'agit de la décomposition primaire (possiblement redondante) de $\overline{I_n}$.
\end{enumerate}
\end{definition}
L'existence en nombre fini de ces valuations et leur unicité 
est l'objet de \cite[Theorem 2.1]{swanson=rees-val} 
et \cite[Theorems 10.1.6 and 10.2.2]{huneke-swanson=integral-closure}. 
Elle repose essentiellement sur le Théorème de Mori--Nagata 
\cite[Theorem 33.10]{nagata_local-rings} déjà cité précédemment. 

D'après la construction dans \cite[Section 10.2]{huneke-swanson=integral-closure} 
et \cite[Alternative construction p.7-8]{swanson=rees-val}, 
les valuations de Rees de $I$ 
sont les valuations associées aux composantes exceptionnelles de $Y\mapto W_0=\Spec(R) $, 
l'éclatement normalisé de $I$.

\begin{propo}\label{propo=Rees}
Avec les hypothèses et notations ci-dessus, 
les diviseurs dicritiques sont les valuations de Rees 
de l'idéal $I:=(f,g)$ où $z=f/g$ et $f,g\in R$ premiers entre eux.
\end{propo}

\begin{demo}
On pourrait se contenter de citer \cite[ex.14.18, p.281]{huneke-swanson=integral-closure}. 
Dans le but d'être le plus complet et le plus géométrique possible, 
nous proposons une démonstration. 
Soit $\gamma: Y\mapto W_0=\Spec(R)$ l'éclatement normalisé de $I$. 
Dans $Y$, l'idéal $I$ est principalisé, 
on a un morphisme $\phi: Y\mapto \PP^1$, $\phi(x):=(f(x):g(x))$.
De même en reprenant la suite d'éclatements~(\ref{eq=1}), on a un morphisme $\psi: W \mapto \PP^1$.
Par propriété universelle de l'éclatement et de la normalisation, 
on a un morphisme $\delta: W \mapto Y$ avec $\psi=\phi \circ \delta$.

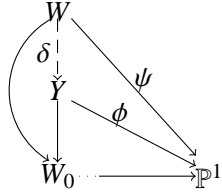
\begin{center}
\begin{pgfpicture}{0.1558cm}{0.4695cm}{3.1558cm}{3.4695cm}%
\pgfsetxvec{\pgfxy(0.9091,0)}
\pgfsetyvec{\pgfxy(0,0.9091)}
\pgfsetstrokecolor{rgb,1:red,0;green,0;blue,0}
\pgfsetlinewidth{0.4pt} 
\pgfputat{\pgfxy(1,1)}{\pgftext{\color{rgb,1:red,0;green,0;blue,0}\small $W_0$}}\pgfstroke
\pgfputat{\pgfxy(1,2.25)}{\pgftext{\color{rgb,1:red,0;green,0;blue,0}\small $Y$}}\pgfstroke
\pgfsetroundcap \pgfsetroundjoin \pgfsetdash{{0pt}{3pt}}{0pt}
\pgfxyline(1.2,1)(1.6,1)
\pgfsetbuttcap \pgfsetdash{}{0pt}
\pgfsetarrows{-to}
\pgfxyline(1.6,1)(3,1)
\pgfsetarrows{-}
\pgfsetdash{{5pt}{3pt}}{0pt}
\pgfsetarrows{-to}
\pgfxyline(1,3.3)(1,2.4)
\pgfsetarrows{-}
\pgfsetdash{}{0pt}
\pgfsetarrows{-to}
\pgfxyline(1.2,2.1)(3,1.15)
\pgfsetarrows{-}
\pgfputat{\pgfxy(3.2,1)}{\pgftext{\color{rgb,1:red,0;green,0;blue,0}\small $\mathbb{P}^1$}}\pgfstroke
\pgfsetarrows{-to}
\pgfxyline(1,2.1)(1,1.2)
\pgfsetarrows{-}
\pgfputat{\pgfxy(1,3.4)}{\pgftext{\color{rgb,1:red,0;green,0;blue,0}\small $W$}}\pgfstroke
\pgfsetarrows{-to}
\pgfxyline(1.2,3.3)(3,1.3)
\pgfsetarrows{-}
\pgfputat{\pgfxy(1.9,1.9)}{\pgftext{\color{rgb,1:red,0;green,0;blue,0}\small $\phi$}}\pgfstroke
\pgfputat{\pgfxy(2.25,2.4)}{\pgftext{\color{rgb,1:red,0;green,0;blue,0}\small $\psi$}}\pgfstroke
\pgfputat{\pgfxy(0.8,2.8)}{\pgftext{\color{rgb,1:red,0;green,0;blue,0}\small $\delta$}}\pgfstroke
\pgfsetarrows{-to}
\pgfpathmoveto{\pgfxy(0.9137,3.297)}\pgfpatharc{119.2488}{238.6713}{1.0909cm}
\pgfstroke
\pgfsetarrows{-}
\end{pgfpicture}%

\end{center}

Donc les dicritiques sont parmi les diviseurs exceptionnels de $Y\mapto W_0=\Spec(R)$. 
Il n'y a plus qu'à montrer que $\phi$ ne contracte aucun de ces diviseurs exceptionnels.

Supposons le contraire, soit $E$ une composante contractée par $\phi$. 
Alors il existe une factorisation $\varepsilon: Y\mapto Y_0 \mapto W_0$, 
avec $Y_0$ normale, où $Y\mapto Y_0$ ne contracte que $E$ 
\cite[p. 238 Correspondance between complete ideals and exceptional curves]{lipman=singrat}. 
Par la propriété universelle de l'éclatement 
$\phi': Y_0\maptodots \PP^1$ a un point fondamental $Q_0$ 
qui est nécessairement l'image de $E$ puisque $Y\mapto \PP^1$ est partout définie. 

On note $\phi (E)=:P \in \PP^1$. 
En termes d'anneaux, le morphisme $\phi$ donne
$$\O_{\PP^1,P}\subseteq\O_{Y_E} \subseteq k(W_0).$$

Comme $\phi'$ est définie hors de $Q_0$, on a
$$O_{\PP^1,P} \subseteq \bigcap_{\height \P_0=1}{O_{Y_0,\P_0}}$$
où l'intersection est sur tous les idéaux premiers $\P_0$ de hauteur~$1$ de $O_{Y_0,Q_0}$. 
En effet, un tel $\P_0$ est l'image par le morphisme \emph{propre} (et birationnel) 
$Y\mapto Y_0$ d'un idéal premier $\P$ d'un certain $O_{Y,Q}$, avec $Q\in E$; 
et on a bien sûr
$O_{P^1,P} \subseteq O_{Y,Q} \subseteq O_{Y,\P}=O_{Y_0,\P_0}$.

Comme $Y_0$ est normale, par le lemme principal des fonctions holomorphes de Zariski
$$O_{Y_0,Q_0}=\bigcap_{\height\P_0=1}{O_{Y_0,\P_0}}.$$
Cela signifie que $\O_{P^1,P}\subseteq\O_{Y_0,Q_0}$, 
\cad/ que l'application $\phi'$ est définie en $Q_0$~: une contradiction.

\end{demo}

Le lecteur remarquera que les arguments ci-dessus 
donnent une troisième preuve de l'existence et de la finitude des dicritiques.


\section{Cas polynomial}
Revenons au cas historique, \cad/ à l'étude des pinceaux de courbes planes.

\begin{definition}\label{defi=dicrit1}(Première définition.)
Soit $k$ un corps, soit $\PP_k^2$ le plan projectif sur $k$, 
et soient deux polynômes $F,G \in k[X,Y,Z]$ homogènes de même degré $d>0$, 
 \emph{premiers entre eux}. 
Le pinceau $\L:=\lambda F + \mu G$ a des \pbs/, 
mais, quitte à faire une composition d'éclatements de points fermés 
$$\Pi:W \mapto \PP_k^2,$$
$\Pi^*(\L)$ définit un morphisme projectif $p:W \mapto \PP_k^1$.
Restreint aux composantes exceptionnelles de $\Pi$, 
on a que $p$ est soit constant, soit surjectif.
Les diviseurs dicritiques de $\L$ 
sont les anneaux locaux des points génériques des composantes exceptionnelles 
où $p$ est surjectif. 
Par abus simplificatif, les composantes exceptionnelles où $p$ est surjectif 
seront appelées aussi diviseurs dicritiques (au sens géométrique).
\end{definition}

Ces diviseurs sont appelés «horizontaux» 
dans \cite[Definition 4]{campillo-piltant-reguera=cones-curves2}. 
Cette définition semble dépendre du choix de $\Pi$. 
Il n'en est rien. 

\begin{propo}\label{propo=dicritPi}
Soient $\Pi_i:W_i \mapto \PP_k^2$ ($i=1,2$) 
deux compositions d'éclatements de points fermés 
telles que $\Pi_i^*(\L)$ définit un morphisme projectif $p_i:W_i \mapto \PP_k^1$. 
Alors les diviseurs dicritiques sont les mêmes pour $\Pi_1$ et $\Pi_2$.
\end{propo}

\begin{demo}
$\Pi_2$ est l'éclatement d'un idéal $I$ de $\O_{\PP_k^2}$; 
quitte à rajouter des éclatements de points fermés, 
on peut supposer que $\Pi_1^{-1}(I)$ est principal, 
\cad/ qu'il existe un morphisme projectif $\Pi_{1,2}: W_1 \mapto W_2$ 
tel que $\Pi_1= \Pi_2 \circ \Pi_{1,2}$. 
Bien sûr, on a $p_1= p_2 \circ \Pi_{1,2}$. 
Donc $p_1$ est constant sur les diviseurs exceptionnels de $\Pi_1$ 
dont l'image par $\Pi_{1,2}$ est un point fermé. 
Les diviseurs dicritiques (au sens géométrique) de $\Pi_1$ 
sont les transformés stricts des diviseurs dicritiques de $\Pi_2$, 
les anneaux locaux en leurs points génériques sont donc les mêmes.
\end{demo}

\begin{definition}\label{defi=dicrit2}(Deuxième définition.)
Soit $k$ un corps et soient $f,g \in k[X,Y]\setminus k$ deux polynômes non constants. 
On note $z:=\frac{f}{g}\in k(X,Y)$. 
Un diviseur dicritique de $z$ 
est un anneau de $k$-valuation discrète $V$ de corps des fractions $k(X,Y)$ 
et tel que le résidu de $z$ dans $K_v:=V/\m_v$ est transcendant sur $k$.
\end{definition}

Pour tout $x\in \PP_k^2$, \pb/ de $\L=\lambda F+\mu G$ 
où $F,G$ sont des homogénéisés de $f,g$ premiers entre eux de même degré, 
on pose $R_x:=\O_{\PP_k^2,x}$. 
On s'aperçoit que les diviseurs dicritiques de $z$ au sens de~\ref{defi=dicrit2} 
sont les diviseurs dicritiques de $z$ pour tous les $R_x$ au sens de~\ref{defi=dicrit}.

\begin{propo}\label{propo=k[z]valuation}
Soit $V$ un anneau de $k$-valuation discrète de $k(X,Y)$ 
tel que l'extension résiduelle $k\mapto V/\m_v$ est transcendante.
On a équivalence:
$$V \text{ est dicritique pour }z\iff V \text{ est un anneau de $k(z)$-valuation.}$$
\end{propo}

\begin{demo}
Bien sûr, si $V$ est dicritique, 
le résidu de tout élément non nul de $k[z]$ est non nul dans $V/\m_v$, 
donc cet élément est de valuation nulle. La réciproque est claire.
\end{demo}

On passe de la première définition à la deuxième en prenant $z=\frac{F}{G}$, 
et de la deuxième à la première en prenant pour $F,G$ des homogénéisés de $f,g$ de même degré. 
Montrons qu'alors les deux définitions sont équivalentes. 
C'est l'objet de la proposition qui suit.

\begin{propo}\label{propo=equivalence-dicrit}
Un diviseur dicritique pour $\L:=\lambda F + \mu G$ (\ref{defi=dicrit1}) 
est dicritique pour $z=\frac{F}{G}$ (\ref{defi=dicrit2}). 
Réciproquement, un diviseur dicritique pour $z=\frac{f}{g}$ (\ref{defi=dicrit2}) 
est dicritique pour $\L:=\lambda F + \mu G$ (\ref{defi=dicrit1}) 
où $F,G$ sont des homogénéisés de $f,g$ de même degré. 
\end{propo}

\begin{demo}
Montrons l'implication directe. 
Un diviseur dicritique pour $\L:=\lambda F + \mu G$ (\ref{defi=dicrit1}) 
est l'anneau local au point générique $\eta$ d'un diviseur $D$ d'une surface régulière~: 
c'est un anneau de valuation discrète $V$. 
Comme $z=\frac{F}{G}$ est défini sauf en un nombre fini de points fermés de $D$, on a $z\in V$. 
Soit un ouvert affine $U$ contenant $\eta$ et où $z$ est défini, 
$\O_U / I(D)$ est le localisé d'un anneau de polynômes $k'[T]$ 
où $k'$ est une extension algébrique de $k$. 
Le résidu de $z$ dans $\O_U / I(D)$ est non constant~: 
il est transcendant sur $k'$ et donc sur $k$, 
comme $K_v:=V/\m_v$ est le corps de fractions de $\O_U / I(D)$, 
on a que $V$ est dicritique pour $z$ au sens de~\ref{defi=dicrit2}.

Réciproquement, soit $V$ un diviseur dicritique pour $z$. 
Il existe $\Pi:W \mapto \PP_k^2$ composition d'éclate\-ments de points fermés 
telle que $\Pi^*(\L)$ définit un morphisme projectif $p:W \mapto \PP_k^1$. 
On conclut en reprenant l'argument de la fin de la première preuve de~\ref{propo=finitude-dicrit}.
\end{demo}

\begin{propo}\label{propo=existence-dicrit}
Avec les hypothèses et notations de~\ref{defi=dicrit1}, 
chaque \pb/ $x$ de $\L:=\lambda F + \mu G$ est centre d'au moins un diviseur dicritique, 
\cad/ qu'il existe au moins un dicritique $V$ 
tel que $\O_{\PP_k^2,x}\subseteq V$ et $\m_v \cap\O_{\PP_k^2,x}= \m_{\PP_k^2,x}$.
\end{propo}


C'est un corollaire de~\ref{propo=finitude-dicrit}.

\begin{propo}[Abhyankar]\label{propo=existence-dicrit2}
Soit $\CC \subseteq \PP_{k(z)}^2$ la courbe générique de $\L$ \cite[p.~517]{campillo-piltant-reguera=cones-curves2}, 
\cad/ la courbe d'équation $F(U,V,T)-zG(U,V,T)\in k(z)[U,V,T]$. 
Les diviseurs dicritiques pour $z$ 
sont les anneaux locaux des points fermés de la désingularisée $\tC$ 
dominant les points d'intersection de $\CC$ et $G(U,V,T)=0$.
\end{propo}

Dans \cite{abh-luengo_dicrit-div}, 
les auteurs regardent le pinceau $\lambda F(U,V,T) + \mu T^N$ 
où $F$ est l'homogénéisé de $f\in k[X,Y]$, $X=U/T,Y=V/T$, $f= F(U,V,T) /T^N$. 
Le pinceau définit en dehors des \pbs/ 
une application $\PP^2(k) \mapto \PP^1$ par $(u,v,t)\mapto (F(u,v,t):t^N)$; 
par restriction, on a un morphisme 
$\Lambda:~\A^2\mapto \A^1=\Spec(k[a])$, $(x,y)\in \A^2 \mapto f(x,y)$. 
Désignons par $\CC$ la courbe affine sur $k(f)$ 
dont l'anneau est $B_f=k[X,Y] \otimes k(f)=k[X,Y] \otimes_{k[a]} k(a)$, 
\cad/ la fibre de $\Lambda$ au dessus du point générique de $\A^1$: 
par définition, $\CC$ est la courbe générique de ce pinceau. 
L'ensemble des diviseurs premiers de $k(X,Y)$ sur $k(f)$, 
noté $D(L/k(f))$ dans \cite{abh-luengo_dicrit-div}, 
est l'ensemble des points de la surface de Riemann (sur $k(f)$) de $\CC$.
L'assertion \cite[(6.2)]{abh-luengo_dicrit-div} signifie que les diviseurs dicritiques de $f$ 
sont les valuations dominant les points à l'infini de $\CC$. 
Ce qui prouve que ces diviseurs existent et sont en nombre fini.
Nous généralisons ici en prenant un pinceau $\lambda F(U,V,T) + \mu G(U,V,T)$.

\begin{demo}
Le polynôme $F(U,V,T)-zG(U,V,T)\in k(z)[U,V,T]$ est homogène. 
On montre fa\-cile\-ment qu'il est irréductible. 
Il définit dans le plan projectif $\PP_{k(z)}^2$ une courbe irréductible $\CC$. 

Plaçons nous dans la carte affine $T\ne 0$. 
On note $u=U/T$, $v=V/T$, $f(u,v):=F(u,v,1)$ et $g(u,v):=G(u,v,1)$. 
On a donc $z=\frac{f(X,Y)}{g(X,Y)}$. 
On note $\overline u$ et $\overline v$ les résidus de $u,v$ dans l'anneau $ k(z)[u,v]/(f(u,v)-zg(u,v))$. 
On a un morphisme $\phi:k[X,Y] \mapto k(z)[u,v] /(f(u,v)-zg(u,v))$ 
défini par $X \mapsto \overline{u}$ et $Y \mapsto \overline{v}$. 
On montre facilement que ce morphisme est injectif. 
Ce morphisme s'étend aux corps de fractions et il définit un isomorphisme entre les deux corps de fractions
$${\tilde \phi}:k(X,Y) \iso k(z)(\CC)$$
où $k(z)(\CC)$ est le corps de fonctions de $\CC$.

On a ${\tilde \phi}(z)=\frac{{\tilde \phi}(f(X,Y))}{{\tilde \phi}(g(X,Y))}=z$.
Donc $\phi$ s'étend (et ${\tilde \phi}$ se restreint) 
à ${\overline \phi}:k(z)[X,Y] \mapto k(z)[u,v] /(f(u,v)-zg(u,v))$ qui est un isomorphisme. 

D'après un résultat classique de Zariski, 
il y a une une bijection entre les points de la variété de Riemann de $\CC$ 
et les $k(z)$-valuations du corps des fractions de $k(z)[u,v] /(f(u,v)-zg(u,v))$ 
(voir \cite[Theorem 41]{zariski-samuel_commut-alg-II} 
ou la rédaction limpide de \cite[Théorème 7.5]{vaquie=valuations-obergurgl}). 
On est en dimension~$1$, la variété de Riemann de $\CC$ est la désingularisée $\tC$ de $\CC$ 
et la bijection de Zariski est simplement l'application $\tC\owns x\mapsto\O_{\tC ,x}$ 
\cite[Theorem 6.12 (b)]{kunz=intro-alg-curves}. 
Par~\ref{propo=k[z]valuation}, si $V$ est un diviseur dicritique pour $z$, 
alors ${\tilde \phi}(V)$ est un des $\O_{\tC ,x}$.


Soit $V$ un diviseur dicritique dominant un \pb/ $x$ de la carte affine d'anneau $k[X,Y]$. 
On a alors la suite d'inclusions
$$k[X,Y] \maptoinj k[X,Y,z]=k[z][X,Y]\maptoinj V.$$
Or tout élément de $k[z]$ est inversible dans $V$, on peut donc insérer $k(z)[X,Y]$ dans la suite ci-dessus.
\begin{equation}\label{eq=2}
k[X,Y] \mapto k[X,Y,z]=k[z][X,Y]\mapto k(z)[X,Y] \mapto V.
\end{equation}
En utilisant l'isomorphisme ${\overline \phi}:k(z)[X,Y] \mapto k(z)[u,v] /(f(u,v)-zg(u,v))$, 
on voit que tous les diviseurs dicritiques dominant $x$ 
correspondent aux points fermés $y\in\tC$ 
dominant $x$ par l'inclusion $k[X,Y] \mapto k(z)[X,Y]$ de~(\ref{eq=2}) 
et que les $y$ sont les points d'intersection de $\tC$ et $G(U,V,T)=0$ dominant~$x$.
\end{demo}

%
%
%
%
%

%
%

Nous reprenons les notations de~\ref{defi=dicrit1}. 
Soient $\L=\angle{F,G}$ un pinceau, 
et $\Pi : W \mapto \PP^2_k$ éliminant les \pbs/ de $\L: \PP^2_k\maptodots \PP^1$, 
avec $W$ régulière. 
Soient ${C} \subseteq \PP_{k}^2$ la courbe d'équation $G=0$ 
et $\Cred$ la courbe réduite correspondante \cite[II.3, p.82]{hartshorne=geometry}. 
On note $p: W \mapto \PP^1$, $C'$ la transformée stricte de $\Cred$ dans $W$ et $O:=p(C')$. 

Voici une nouvelle généralisation du théorème d'Abhyankar--Luengo 
\cite[Theorems (1.1), (7.1), (7.2), (7.3)]{abh-luengo_dicrit-div}.

\begin{propo}\label{propo=Abhyankar-Luengo} 
On suppose $\Cred$ lisse en les \pbs/ de $\L$. 
Pour toute composante dicritique $D \subseteq W$ de $\L$, 
on a que $D\cap p^{-1}(O)$ se réduit à un point fermé $P$. 
Ainsi $z:=\frac{F}{G}$ peut être défini sur $D\setminus\{P\}$ 
qui est une droite affine dont l'anneau de fonctions est une algèbre de polynômes~: 
$z$ est résiduellement un polynôme au sens de~\cite[Theorem (7.1)]{abh-luengo_dicrit-div}.
\end{propo}

\begin{demo}
On suppose que:
\begin{enumerate}[1)]
\item\label{sup=1}
La fibre $p^{-1}(O)$ est connexe.
\item\label{sup=2}
Soit $\Gamma$ le graphe obtenu comme suit~: 
On prend le graphe dual des composantes de $\Pi^{-1}(\Cred)$ 
et on contracte en un seul point $\Omega$ toutes les composantes irréductibles de $\Cred$. 
Alors $\Gamma$ est un arbre de racine $\Omega$.
\end{enumerate}

Admettons \ref{sup=1} et \ref{sup=2} et prouvons la proposition. 
$\Gamma \setminus \{D\}$ a des composantes connexes $\Gamma_0, \Gamma_1, \ldots , \Gamma_s$, 
et on choisit $\Gamma_0$ pour que $\Omega \in \Gamma_0$. 
Comme $p^{-1}(O)$ est connexe et contient $\Omega$ mais pas $D$, 
on a $p^{-1}(O) \subseteq \Gamma_0$. 
Mais $\Gamma$ est un arbre, $\Gamma_0$ aussi et a pour racine $\Omega$, 
donc $D$ est rattaché à $\Gamma_0$ en exactement un seul point~: 
son prédécesseur que nous notons $E_D$. 
Par ailleurs, $p$ restreinte à $D$ est propre, $p(D)$ non constant, donc $p(D)=\PP^1$. 
Comme $\Gamma$ est un arbre, $\{E_D\}\cap D$ est un point $P$, 
éventuellement l'isomorphe d'un \pb/ par $C'\mapto \Cred$ si $E_D=\Omega$~: 
la courbe $\Cred$ est lisse en les \pbs/ $\PP^2$ 
et est donc isomorphe à son transformé strict $C'$ qui est connexe. 
Ainsi $P=\{E_D\}\cap D\supseteq p^{-1}(O)\cap D\ne \emptyset$: $p^{-1}(O)\cap D$ est réduit à $P$. 
CQFD

Prouvons \ref{sup=2}. 
Soient $\Gamma'_i:=\Pi^{-1}(B_i)$ où $B_i$ ($1\le i \le t$) sont les \pbs/ de $\L=\angle{F,G}$. 
Comme $\Cred$ est lisse, 
$\Gamma'_1, \ldots , \Gamma'_t$ sont des arbres 
et $\Gamma$ est l'arbre obtenu en joignant les racines respectives de $\Gamma'_1, \ldots , \Gamma'_t$ 
(\ie/ les composantes exceptionnelles qui intersectent $C'$) à $\Omega$.

Pour \ref{sup=1}, par la factorisation de Stein \cite[III. Corollary 11.5]{hartshorne=geometry}, 
$p$ se factorise par $p_1: W \mapto Y$ et $n: Y \mapto \PP^1$ 
où $n$ est un morphisme fini, $Y$ normale et $p_1$ est projective et a toutes ses fibres connexes. 
Soit $\Lambda$ une courbe irréductible de $\PP^2$ qui ne passe pas par les \pbs/ de $\L$. 
Par Bézout, $\Lambda$ intersecte toutes les courbes de $\L$. 
Soit $\Lambda'$ sa transformée stricte dans $W$.
On a $p_1(\Lambda')=Y$ parce que $\Lambda$ 
n'est contenue dans aucune courbe de $\L$ et que $p_1$ est fermée. 
On en déduit que $\Lambda'$ 
rencontre toutes les composantes connexes de $p^{-1}(O)=p_1^{-1}(n^{-1}(O))$ 
(elles sont en bijection avec les points de $n^{-1}(O)$). 
Comme $\Lambda$ ne passe pas par les \pbs/ de $\L$, 
les points de $\Lambda'\cap p^{-1}(O)$ ne sont pas exceptionnels pour $W\mapto \PP^2$. 
Cela entraîne que chaque composante connexe de $p^{-1}(O)$ 
contient au moins une composante irréductible non exceptionnelle, 
\cad/ le transformé strict d'une composante irréductible de $\Cred$. 
Par hypothèse, on n'a fait des éclatements qu'au dessus de points lisses de $\Cred$: 
le transformé strict de $\Cred$ est connexe. 
Donc $p^{-1}(O)$ a une seule composante connexe. 
\end{demo}

On extrait de la preuve le résultat géométrique suivant 
qui généralise le théorème de connexité de \cite{le-weber_jacobian-dim2} p.~377.

\begin{coro}\label{coro=Abhyankar-Luengo-geom}
Avec les hypothèses et notations de~\ref{propo=Abhyankar-Luengo}, la fibre $p^{-1}(O)$ est connexe.
\end{coro}

\end{document}